\documentclass{article}

\def\ss{\subseteq}

\def\zetabar{\overline{\zeta}}

\def\dbar{\overline{\partial}}

 \def\HollowBox #1#2{{\dimen0=#1 \advance\dimen0 by -#2
       \dimen1=#1 \advance\dimen1 by #2
        \vrule height #1 depth #2 width #2
        \vrule height 0pt depth #2 width #1
        \llap{\vrule height #1 depth -\dimen0 width \dimen1}%
       \hskip -#2
       \vrule height #1 depth #2 width #2}}
 
\font\teneufm=eufm10
\font\seveneufm=eufm7
\font\fiveeufm=eufm5
\newfam\eufmfam
\textfont\eufmfam=\teneufm
\scriptfont\eufmfam=\seveneufm
\scriptscriptfont\eufmfam=\fiveeufm

\newfam\msbfam
\font\tenmsb=msbm10  \textfont\msbfam=\tenmsb
\font\sevenmsb=msbm7  \scriptfont\msbfam=\sevenmsb
\font\fivemsb=msbm5    \scriptscriptfont\msbfam=\fivemsb
\def\Bbb{\fam\msbfam \tenmsb}

\def\CC{{\Bbb B}}
\def\CC{{\Bbb C}}

\newfam\msbbfam
\font\tenmsbb=msbm10  scaled \magstep2 \textfont\msbbfam=\tenmsbb
\font\sevenmsbb=msbm7 scaled \magstep2 \scriptfont\msbbfam=\sevenmsbb
\font\fivemsbb=msbm5  scaled \magstep2 \scriptscriptfont\msbbfam=\fivemsbb
\def\Bbbb{\fam\msbbfam \tenmsbb}

\def\CCC{{\Bbbb C}}

\newtheorem{theorem}{Theorem}[section]

\makeindex

\usepackage{graphicx}

\usepackage{amsmath}

\begin{document}

\begin{center}
\huge \bf
A Direct Connection Between the Bergman and Szeg\H{o} Kernels\footnote{{\bf MR Classification
Numbers:}  32A25}
\end{center}
\vspace*{.12in}

\begin{center}
\Large Steven G. Krantz
\end{center}
\vspace*{.15in}

\begin{center}
\today
\end{center}
\vspace*{.2in}

\begin{quotation}
{\bf Abstract:} \sl
We use Stokes's theorem to establish an explicit and concrete connection between the
Bergman and Szeg\H{o} projections on the disc, the ball, and on strongly
pseudoconvex domains.  
\end{quotation}
\vspace*{.25in}

\section{Introduction}

Two of the most classical and well established reproducing formulas in 
complex analysis are those of S. Bergman and G. Szeg\H{o}.  The first
of these is a formula for the Bergman space, and the associated integral
lives on the interior of the domain in question.  The latter of these is 
a formula for the Hardy space, and the associated integral lives on the boundary
of the domain.  For formal reasons, the Bergman integral gives rise to a projection
from $L^2(\Omega)$ to $A^2(\Omega)$ (the Bergman space); likewise, the Szeg\H{o} 
integral gives rise to a projection from $L^2(\partial \Omega)$ to $H^2(\Omega)$ (the
Hardy space).  

Since both of the artifacts in question here are canonical, it is natural to suspect
that there is some relationship between the two integral formulas.  After all, they
both reproduce functions that are continuous on the closure of the domain and 
holomorphic on the interior.  In the present paper we establish such a connection---very
explicitly---on a variety of domains in $\CC^1$ and $\CC^n$.  This is done by way
of a moderately subtle calculation using Stokes's theorem.  The calculation itself
has some intrinsic interest, but the main point is the equality of the canonical
integrals and the associated projections.

\section{The Case of the Disc}

Let $D$ be the unit disc in $\CC$.  In this context, the Szeg\H{o} kernel is
$$
S(z, \zeta) = \frac{1}{2\pi} \cdot \frac{1}{1 - z \cdot \overline{\zeta}} 
$$
and the Bergman kernel is 
$$
K(z,\zeta) = \frac{1}{\pi} \cdot \frac{1}{(1 - z \cdot \overline{\zeta})^2} \, .
$$

Take $f$ to be real analytic on a neighborhood of $\overline{D}$.  Now we can calculate
\begin{eqnarray*}
\frac{1}{2\pi} \int_{\partial D} f(\zeta) S(z, \zeta) d\sigma(\zeta) & = &
    \frac{1}{2\pi} \int_{\partial D} f(\zeta) \cdot \frac{1}{1 - z \overline{\zeta}} \, \frac{\left [ \zetabar d\zeta - \zeta d\zetabar \right ]}{2i}  \\
    & = & \frac{1}{4\pi i} \int_{\partial D} \frac{f(\zeta) \zetabar}{1 - z \cdot \zetabar} \, d\zeta 
              - \frac{1}{4\pi i} \int_{\partial D} \frac{f(\zeta)\zeta}{1 - z\cdot \zetabar} \, d\zetabar \\
    & \stackrel{\rm (Stokes)}{=} & \frac{1}{4\pi i} \mathop{\int \!\! \int}_D \frac{f(\zeta)}{1 - z\cdot \zetabar} \, d\zetabar \wedge d\zeta 
         + \frac{1}{4\pi i} \mathop{\int \! \! \int}_D \frac{f(\zeta) \zetabar z}{(1 - z \cdot \zetabar)^2} \, d\zetabar \wedge d\zeta   \\
    && \qquad - \frac{1}{4\pi i} \mathop{\int \! \! \int}_D \frac{\partial (f \cdot \zeta)/\partial \zeta}{1 - z \cdot \zetabar} \, d\zeta \wedge d\zetabar 
    + \frac{1}{4\pi i} \mathop{\int \!\! \int}_D \frac{\partial f/\partial \zetabar \cdot \zetabar}{1 - z \cdot \zetabar} \, d\zetabar \wedge d\zeta \\
    & = & \frac{1}{4\pi i} \mathop{\int \!\! \int}_D \frac{f(\zeta)}{(1 - z \cdot \zetabar)^2} \, d\zetabar \wedge d\zeta 
         - \frac{1}{4\pi i} \mathop{\int \!\! \int}_D \frac{f(\zeta)z\zetabar}{(1 - z \cdot \zetabar)^2} \, d\zetabar \wedge d\zeta   \\
    && \qquad + \frac{1}{4\pi i} \mathop{\int \!\! \int}_D \frac{f(\zeta)\zetabar z}{(1 - z \cdot \zetabar)^2} \, d\zetabar \wedge d\zeta
	 - \frac{1}{4\pi i} \mathop{\int \!\! \int}_D \frac{\partial (f \cdot\zeta)/\partial \zeta}{1 - z \cdot \zetabar} \, d\zeta \wedge d\zetabar \\
    && \qquad + \frac{1}{4\pi i} \mathop{\int \!\! \int}_D \frac{\partial f/\partial \zetabar \cdot \zetabar}{1 - z \cdot \zetabar} \, d\zetabar \wedge d\zeta \\
    & = & \frac{1}{4\pi i} \mathop{\int \!\! \int}_D \frac{f(\zeta)}{(1 - z \cdot \zetabar)^2} \, d\zetabar \wedge d\zeta 
	 - \frac{1}{4\pi i} \mathop{\int \!\! \int}_D \frac{\partial (f \cdot\zeta)/\partial \zeta}{1 - z \cdot \zetabar} \, d\zeta \wedge d\zetabar \\
    && \qquad + \frac{1}{4\pi i} \mathop{\int \!\! \int}_D \frac{\partial f/\partial \zetabar \cdot \zetabar}{1 - z \cdot \zetabar} \, d\zetabar \wedge d\zeta \\
    & = & \frac{1}{4\pi i} \mathop{\int \!\! \int}_D \frac{f(\zeta)}{(1 - z \cdot \zetabar)^2} \, d\zetabar \wedge d\zeta 
	 - \frac{1}{4\pi i} \mathop{\int \!\! \int}_D \frac{\partial f/\partial \zeta \cdot \zeta}{1 - z \cdot \zetabar} \, d\zeta \wedge d\zetabar  \\
    && \qquad - \frac{1}{4\pi i} \mathop{\int \!\! \int}_D \frac{f(\zeta)}{1 - z \cdot \zetabar} \, d\zeta \wedge d\zetabar 
     + \frac{1}{4\pi i} \mathop{\int \!\! \int}_D \frac{\partial f/\partial \zetabar \cdot \zetabar}{1 - z \cdot \zetabar} \, d\zetabar \wedge d\zeta	\\
\end{eqnarray*}

\begin{eqnarray*}
    & = & \frac{1}{4\pi i} \mathop{\int \!\! \int}_D \frac{f(\zeta)}{(1 - z \cdot \zetabar)^2} \, d\zetabar \wedge d\zeta   
	 - \frac{1}{4\pi i} \mathop{\int \!\! \int}_D \frac{\partial f/\partial \zeta \cdot\zeta}{1 - z \cdot \zetabar} \, d\zeta \wedge d\zetabar \\
    && \qquad - \frac{1}{4\pi i} \mathop{\int \!\! \int}_D \frac{f(\zeta)}{(1 - z \cdot \zetabar)^2} \, d\zeta \wedge d\zetabar 
	 + \frac{1}{4\pi i} \mathop{\int \!\! \int}_D \frac{f(\zeta)z \zetabar}{(1 - z \cdot \zetabar)^2} \, d\zeta \wedge d\zetabar \\
    && \qquad + \frac{1}{4\pi i} \mathop{\int \!\! \int}_D \frac{\partial f/\partial \zetabar \cdot \zetabar}{1 - z \cdot \zetabar} \, d\zetabar \wedge d\zeta	\\
    & = & \frac{1}{2\pi i}\mathop{\int \!\! \int}_D \frac{f(\zeta)}{(1 - z \cdot \zetabar)^2} \, d\zetabar \wedge d\zeta 
	 - \frac{1}{4\pi i} \mathop{\int \!\! \int}_D \frac{\partial f/\partial \zeta \cdot \zeta}{1 - z \cdot \zetabar} \, d\zeta \wedge d\zetabar \\
    && \qquad + \frac{1}{4\pi i} \mathop{\int \!\! \int}_D \frac{f(\zeta)z\zetabar}{(1 - z \cdot \zetabar)^2} \, d\zeta \wedge d\zetabar 
     + \frac{1}{4\pi i} \mathop{\int \!\! \int}_D \frac{\partial f/\partial \zetabar \cdot \zetabar}{1 - z \cdot \zetabar} \, d\zetabar \wedge d\zeta \\
    & = & A - B + C + D \, .
\end{eqnarray*}

Certainly $A = \int_D f(\zeta) K(z, \zeta) \, dA(\zeta)$, where $K$ is the Bergman kernel of the disc. 
So this is the Bergman projection.  Now we claim that $- B
+ C + D \equiv 0$. If we can establish that assertion, then we
will have seen directly, by way of Stokes's theorem, that the
Szeg\"{o} projection equals the Bergman projection (at least for functions
real analytic on the closure).

First assume that $f$ is holomorphic.
We establish the claim by verifying it for $f(\zeta) = \zeta^k$, each $k = 0, 1, 2, \dots$.
Indeed, in this case (expanding the kernel in a Neumann series and discarding terms that
obviously integrate to zero by parity)
$$
B = \frac{1}{4\pi i} \mathop{\int \!\! \int}_D k z^k |\zeta|^{2k} \, d\zetabar \wedge d\zeta \, .
$$
And a similar calculation shows that
$$
C = \frac{1}{4\pi i} \mathop{\int \!\! \int}_D k z^k |\zeta|^{2k} \, d\zetabar \wedge d\zeta \, .
$$
And $D = 0$ because $\partial f/\partial \zetabar \equiv 0$.  Thus $-B + C + D = 0$ as desired.

For any monomial containing some positive power of $\overline{\zeta}$, it is easy
to see by parity (again using the Neumann series for the kernel) that the integrals
$B, C, D$ are  equal to 0.
Summing, we see that we have proved our result for any function $f$ that is real analytic
on a neighborhood of $\overline{D}$.  But standard measure theory, together with
the Weierstrass approximation theorem, enable us to pass from these functions to,
for example, functions that are continuous on $\overline{D}$.

Thus we see by our calculation that the full Szeg\H{o} projection is equal to the full
Bergman projection on the disc $D$.

We treat the case of the Bergman and Szeg\H{o} projections on the ball below.

Given Fefferman's asymptotic expansion for the Bergman kernel [FEF], and Boutet de Monvel/Sj\"{o}strand's
asymptotic expansion for the Szeg\H{o} kernel [BOS], one would expect a like calculation (up to a
controllable error term) on a smoothly bounded, strongly pseudoconvex domain.  Unfortunately we
do not know enough about the canonical kernels on domains of finite type to be able to 
predict what will happen there.  We explore the strongly pseudoconvex case below.

In a more recent work, Chen and Fu [CHF] have explored some new comparisons of the Bergman
and Szeg\"{o} kernels.  A sample theorem is this:

\begin{theorem} \sl 
Let $\Omega \ss \CC^n$ be a pseudoconvex domain with $C^2$ boundary.  Then
\begin{enumerate}
\item[{\bf (1)}] For any $0 < a < 1$, there exists a constant $C > 0$ such that
$$
\frac{S(z,z)}{K(z,z)} \leq C \delta(z) |\log \delta(z) |^{n/a} \, .
$$
\item[{\bf (2)}]  If there is a neighborhood $U$ of $\partial \Omega$, a bounded, continuous 
plurisubharmonic function $\varphi$ on $U \cap \Omega$, and a defining function $\rho$ of $\Omega$
satisfying $i \partial \dbar \varphi \geq i \rho^{-1} \partial \dbar \rho$ on $U \cap \Omega$ as currents,
then there exists constants $0 < a < 1$ and $C > 0$ such that
$$
\frac{S(z,z)}{K(z,z)} \geq C \delta(z) |\log \delta(z)|^{-1/a} \, .
$$
\end{enumerate}
\end{theorem}

These authors further show that, for a $C^2$-bounded {\it convex} domain the quotient $S/K$ is comparable
to $\delta$ without any logarithmic factor.

The techniques used in this work are weighted estimates for the $\dbar$ operator (in the spirit
of H\"{o}rmander's work [HOR]) and also an innovative use of the Diederich-Forn\ae ss index
(see [DIF]).   We can say no more about the work here.

We turn next to an examination of the situation on the unit ball $B$ in $\CC^n$.

\section{The Unit Ball in \boldmath $\CCC^n$}

For simplicity we shall in fact restrict attention to complex dimension 2.  In that situation, the
area measure $d\sigma$ on the boundary is given by
$$
d\sigma =  \frac{1}{16} \, \biggl [ \zeta_1 d\zeta_2 \wedge d\zetabar_1 \wedge d\zetabar_2 
    	     - \zeta_2 d\zeta_1 \wedge d\zetabar_1 \wedge d\zetabar_2 
	     + \zetabar_1 d\zetabar_2 \wedge d\zeta_1 \wedge d\zeta_2 
	     - \zetabar_2 d\zetabar_1 \wedge d\zeta_1 \wedge d\zeta_2 \biggr ]  \, .
$$
As a result, we have
\begin{eqnarray*}
\int_{\partial B} f(\zeta) S(z, \zeta) \, d\sigma(\zeta) & = & 
      \frac{1}{2\pi^2} \cdot \frac{1}{16} \mathop{\int \!\!\! \int \!\!\! \int}_{\partial B} f(\zeta) \frac{1}{(1 - z \cdot \zetabar)^2} \, 
                      \biggl [ \zeta_1 d\zeta_2 \wedge d\zetabar_1 \wedge d\zetabar_2 
    	     - \zeta_2 d\zeta_1 \wedge d\zetabar_1 \wedge d\zetabar_2  \\
      && \qquad	+ \zetabar_1 d\zetabar_2 \wedge d\zeta_1 \wedge d\zeta_2 
	     - \zetabar_2 d\zetabar_1 \wedge d\zeta_1 \wedge d\zeta_2 \biggr ] \\
      & = & \frac{1}{32\pi^2} \mathop{\int \!\!\! \int \!\!\! \int \!\!\! \int}_B \frac{\partial f/\partial \zeta_1}{(1 - z \cdot \zetabar)^2} \cdot \zeta_1 \,
                  d\zeta_1 \wedge d\zeta_2 \wedge d\zetabar_1 \wedge d\zetabar_2  \\
      \end{eqnarray*}

\begin{eqnarray*}
&& \qquad  + \frac{1}{32\pi^2} \mathop{\int \!\!\! \int \!\!\! \int \!\!\! \int}_B \frac{f}{(1 - z \cdot \zetabar)^2} \, d\zeta_1 \wedge d\zeta_2 \wedge d\zetabar_1 \wedge d\zetabar_2 \\
      && \qquad -  \frac{1}{32\pi^2} \mathop{\int \!\!\! \int \!\!\! \int \!\!\! \int}_B \frac{\partial f/\partial \zeta_2}{(1 - z \cdot \zetabar)^2} \cdot \zeta_2 \,
                  d\zeta_2 \wedge d\zeta_1 \wedge d\zetabar_1 \wedge d\zetabar_2 \\
      && \qquad  - \frac{1}{32\pi^2} \mathop{\int \!\!\! \int \!\!\! \int \!\!\! \int}_B \frac{f}{(1 - z \cdot \zetabar)^2} \, d\zeta_2 \wedge d\zeta_1 \wedge d\zetabar_1 \wedge d\zetabar_2 \\
      && \qquad + \frac{1}{32\pi^2} \mathop{\int \!\!\! \int \!\!\! \int \!\!\! \int}_B \frac{\partial f/\partial \zetabar_1}{(1 - z \cdot \zetabar)^2} \cdot \zetabar_1 \,
                d\zetabar_1 \wedge d\zetabar_2 \wedge d\zeta_1 \wedge d\zeta_2 \\
      && \qquad + \frac{1}{32\pi^2} \mathop{\int \!\!\! \int \!\!\! \int \!\!\! \int}_B \frac{f}{(1 - z \cdot \zetabar)^2} \, d\zetabar_1 \wedge d\zetabar_2 \wedge d\zeta_1 \wedge d\zeta_2 \\
      && \qquad + \frac{1}{32\pi^2} \mathop{\int \!\!\! \int \!\!\! \int \!\!\! \int}_B \frac{2f \cdot \zetabar_1 z_1}{(1 - z \cdot \zetabar)^3} \, d\zetabar_1 \wedge d\zetabar_2 \wedge d\zeta_1 \wedge d\zeta_2  \\
      && \qquad - \frac{1}{32\pi^2} \mathop{\int \!\!\! \int \!\!\! \int \!\!\! \int}_B \frac{\partial f/\partial \zetabar_2}{(1 - z \cdot \zetabar)^2} \cdot \zetabar_2 \,
           d\zetabar_2 \wedge d\zetabar_1 \wedge d\zeta_1 \wedge d\zeta_2 \\
      && \qquad  - \frac{1}{32\pi^2} \mathop{\int \!\!\! \int \!\!\! \int \!\!\! \int}_B \frac{f}{(1 - z \cdot \zetabar)^2} \, d\zetabar_2 \wedge d\zetabar_1 \wedge d\zeta_1 \wedge d\zeta_2 \\
      && \qquad  - \frac{1}{32\pi^2} \mathop{\int \!\!\! \int \!\!\! \int \!\!\! \int}_B \frac{2f \cdot \zetabar_2 z_2}{(1 - z \cdot \zetabar)^3} \, d\zetabar_2 \wedge d\zetabar_1 \wedge d\zeta_1 \wedge d\zeta_2 \, . \\
\end{eqnarray*}

\noindent Now we may group together like terms to obtain

\begin{eqnarray*}
& = & - \frac{1}{8\pi^2} \mathop{\int \!\!\! \int \!\!\! \int \!\!\! \int}_B \frac{f(\zeta)}{(1 - z \cdot \zetabar)^3} \, d\zetabar_1 \wedge d\zeta_1 \wedge d\zetabar_2 \wedge d\zeta_2 \\
   && \qquad + \frac{3}{16\pi^2} \mathop{\int \!\!\! \int \!\!\! \int \!\!\! \int}_B \frac{f(\zeta) \cdot (z \cdot \zetabar)}{(1 - z \cdot \zetabar)^3} \, d\zetabar_1 \wedge d\zeta_1 \wedge d\zetabar_2 \wedge d\zeta_2 \\
   && \qquad + \frac{1}{32\pi^2} \mathop{\int \!\!\! \int \!\!\! \int \!\!\! \int}_B \frac{\partial f/\partial \zeta_1}{(1 - z \cdot \zetabar)^2} \cdot \zeta_1 \,
                  d\zeta_1 \wedge d\zeta_2 \wedge d\zetabar_1 \wedge d\zetabar_2 \\
   && \qquad -  \frac{1}{32\pi^2} \mathop{\int \!\!\! \int \!\!\! \int \!\!\! \int}_B \frac{\partial f/\partial \zeta_2}{(1 - z \cdot \zetabar)^2} \cdot \zeta_2 \,
                  d\zeta_2 \wedge d\zeta_1 \wedge d\zetabar_1 \wedge d\zetabar_2 \\
\end{eqnarray*}

\begin{eqnarray*}
   && \qquad + \frac{1}{32\pi^2} \mathop{\int \!\!\! \int \!\!\! \int \!\!\! \int}_B \frac{\partial f/\partial \zetabar_1}{(1 - z \cdot \zetabar)^2} \cdot \zetabar_1 \,
                d\zetabar_1 \wedge d\zetabar_2 \wedge d\zeta_1 \wedge d\zeta_2 \\
   && \qquad - \frac{1}{32\pi^2} \mathop{\int \!\!\! \int \!\!\! \int \!\!\! \int}_B \frac{\partial f/\partial \zetabar_2}{(1 - z \cdot \zetabar)^2} \cdot \zetabar_2 \,
           d\zetabar_2 \wedge d\zetabar_1 \wedge d\zeta_1 \wedge d\zeta_2 \\
   & = & - A + B + C - D + E - F \, .
\end{eqnarray*}

Now $-A$ is just the usual Bergman integral on the ball $B$ in $\CC^2$.   And we can argue, just as on the disc, that
the other terms cancel out (or are zero outright, just by parity).  We have verified that
the Szeg\H{o} projection integral equals the Bergman projection integral on the unit
ball $B \ss \CC^n$.

\section{Strongly Pseudoconvex Domains}

We again, for simplicity, restrict attention to $\CC^2$. In the seminal
paper [FEF], Fefferman shows that, near a strongly pseudoconvex boundary
point, the Bergman kernel may be written (in suitable local coordinates)
as
$$
\frac{2}{\pi^2} \cdot \frac{1}{(1 - z \cdot \zetabar)^3} + {\cal E}(z, \zeta) \, ,
$$
where ${\cal E}$ is an error term of strictly lower order (in the sense of pseudodifferential
operators) than the Bergman kernel.

In the important paper [BOS], Boutet de Monvel and Sj\"{o}strand show that, near a strongly 
pseudoconvex boundary point, the Szeg\H{o} kernel may be written (in suitable local coordinates) as
$$
\frac{1}{2\pi^2} \cdot \frac{1}{(1 - z \cdot \zetabar)^2} + {\cal F}(z, \zeta) \, ,
$$
where ${\cal F}$ is an error term of strictly lower order (in the sense of pseudodifferential
or Fourier integral operators) than the Szeg\H{o} kernel.

We now take advantage of these two asymptotic expansions to say something about the relationship
between the Bergman and Szeg\H{o} projections on a smoothly bounded strongly pseudoconvex
domain.

Now fix a smoothly bounded, strongly pseudoconvex domain $\Omega$ with
defining function $\rho$ (see [KRA1] for this notion). Let $U$ be a
tubular neighborhood of $\partial \Omega$ and let $V$ be a relatively
compact subdomain of $U$ that is also a tubular neighborhood of $\partial \Omega$.  Let $\varphi_j$
be a partition of unity that is supported in $U$ and sums to be identically 1
on $V$.  We assume that each $\varphi_j$ has support so small that both
the Fefferman and Boutet de Monvel/Sj\"{o}strand expansions are valid
on the support of $\varphi_j$.  Then we write
\begin{eqnarray*}
 \int_{\partial \Omega} f(\zeta) S(z, \zeta) \, d\sigma(\zeta) & = & 
      \mathop{\int \!\!\! \int \!\!\! \int}_{\partial \Omega} f(\zeta) S(z, \zeta) \omega(\zeta) \\
   & = & \sum_j \mathop{\int \!\!\! \int \!\!\! \int}_{\partial \Omega} \varphi_j(\zeta) f(\zeta) S(z, \zeta) \omega(\zeta) \, ,
\end{eqnarray*}
where $\omega$ is the differential form that is equivalent to area measure on the boundary.
And now, using Boutet de Monvel/Sj\"{o}strand, and using the notable lemma of Fefferman [FEF] that
says that a strongly pseudoconvex boundary point is the ball up to fourth order, one can write each
term of this last sum as 
\begin{eqnarray*}
\lefteqn{\frac{1}{2\pi} \mathop{\int \!\!\! \int \!\!\! \int}_{\partial B} \widetilde{\varphi}_j(\zeta) f(\zeta) \cdot \frac{1}{(1 - z \cdot \zetabar)^2} 
    \biggl [ \zeta_1 d\zeta_2 \wedge d\zetabar_1 \wedge d\zetabar_2}  \\ 
    && \qquad - \zeta_2 d\zeta_1 \wedge d\zetabar_1 \wedge d\zetabar_2  
      	+ \zetabar_1 d\zetabar_2 \wedge d\zeta_1 \wedge d\zeta_2 
	     - \zetabar_2 d\zetabar_1 \wedge d\zeta_1 \wedge d\zeta_2 \biggr ]  + {\cal G} \, ,
\end{eqnarray*}
where the error term ${\cal G}$ arises from approximating $\partial \Omega$ by $\partial B$,
from approximating the Szeg\H{o} kernel $S$ by
the kernel for the ball, by applying a change of variable to $\varphi_j$, 
and also by approximating $\omega$ by the differential
form that we used on the ball.

Now we may carry out the calculations using Stokes's theorem just as in the last section to
finally arrive at the assertion that the last integral equals
$$
\mathop{\int \!\!\! \int \!\!\! \int \!\!\! \int}_B \widetilde{\widetilde{\varphi}}_j(\zeta) \frac{f(\zeta)}{(1 - z \cdot \zetabar)^3} \, dV +{\cal H} \, .
$$
We cannot make the error term ${\cal H}$ disappear this time, but it is smoothly bounded hence negligeble.  Finally, we can use
the Fefferman asymptotic expansion to relate this last integral to the Bergman projection integral on the
strongly pseudoconvex domain $\Omega$.

In summary, we have used Stokes's theorem to relate the Szeg\H{o} projection integral on a smoothly
bounded, strongly pseudoconvex domain to the Bergman projection integral on that domain.  In this
context, we do not get a literal equality.  Instead we get an equality up to a controllable error term.

\section{Concluding Remarks}

Certainly one of the fundamental problems of the function theory of several complex variables is
to understand the canonical kernels in as much detail as possible.  This paper is a contribution
to that program.   In future papers we hope to explore the finite type case in $\CC^n$ and other
more general domains as well.

\newpage

\noindent {\Large \sc References}
\bigskip  \\

\begin{enumerate} 

\item[{\bf [BOS]}]  L. Boutet de Monvel and J. Sj\"{o}strand, 
Sur la singularit\'{e} des noyaux de Bergman et Szeg\"{o},
{\it Soc.\ Mat.\ de France Asterisque} 34-35(1976), 123--164.

\item[{\bf [CHF]}] B.-Y. Chen and S. Fu, Comparison of the
Bergman and Szeg\"{o} kernels, {\it Advances in Math.}
228(2011), 2366--2384.

\item[{\bf [FEF]}] C. Fefferman, The Bergman kernel and
biholomorphic mappings of pseudoconvex domains, {\it Invent. Math.}
26(1974), 1--65.

\item[{\bf [KRA1]}]  S. G. Krantz, {\it Function Theory of Several Complex
Variables}, 2nd ed., American Mathematical Society, Providence, RI, 2001.

\end{enumerate}
\vspace*{.17in}

\begin{quote}
Department of Mathematics \\
Washington University in St.\ Louis  \\
St.\ Louis, Missouri 63130 \ \ U.S.A.  \\
{\tt sk@math.wustl.edu}
\end{quote}

\end{document}